\documentclass[12pt,leqno]{article}
\usepackage[francais]{babel}
\usepackage{amsmath,amsfonts,amssymb,amsthm}
\newtheorem{montheo}{Th\'eor\`eme}

\newcommand\Z{\mathbb{Z}}

\newcommand\N{\mathbb{N}}
\newcommand\R{\mathbb{R}}
\newcommand\rg{\rightarrow}

\begin{document}

\title{Mailles et ensembles de Sidon}
\author{Jean--Pierre Kahane}
\date{}
\maketitle

Le terme  d'ensemble de Sidon est apparu en 1957, en liaison avec une \'etude sur les fonctions moyenne--p\'eriodiques born\'ees \cite{kahane}. Une propri\'et\'e en \'etait signal\'ee comme ``condition de maille'',  et je me suis souvent demand\'e si cette condition \'etait am\'eliorable, ou si elle \'etait n\'ecessaire et suffisante. Le pr\'esent article r\'epond n\'egativement \`a ces deux questions.

La partie 1 contient les d\'efinitions et les principaux \'enonc\'es. La partie~2 donne une construction d'ensembles quasi--ind\'ependants qui \'etablit que la condition de maille est inam\'eliorable. Les parties~3, 4 et 5 montrent que, m\^eme  consid\'erablement renforc\'ee, elle est loin de garantir qu'un ensemble est de Sidon ; ces parties font appel \`a l'outil de s\'election al\'eatoire introduit par Katznelson et Malliavin en 1966 \cite{katmal}, et consid\'erablement d\'evelopp\'e par Bourgain dans sa th\'eorie des d\'efinitions \'equivalentes des ensembles de Sidon \cite{bourg} \cite{lique} ; elles se r\'ef\`erent pour l'essentiel \`a l'\'etude  des ensembles de Sidon faite par Pisier en 1981 \cite{pisier}. La partie~6 lie condition de maille et ensembles d'analyticit\'e, en s'inspirant de \cite{katmal}. Un appendice d\'etaille les calculs de probabilit\'es utilis\'es dans l'article.

Mon int\'er\^et pour les ensembles de Sidon s'est r\'eveill\'e \`a l'occasion du colloque organis\'e \`a Orsay en janvier 2005 en l'honneur de Myriam D\'echamps. Les travaux de Myriam D\'echamps appartiennent \`a l'histoire des ensembles de Sidon, qu'il s'agisse de contributions originales ou de mises au point \cite{dec2} \cite{dec} \cite{lique}. Lors du colloque j'avais annonc\'e sans en avoir la preuve que la condition de maille n'\'etait pas suffisante pour avoir un ensemble de Sidon. La mise au point a \'et\'e laborieuse et elle a b\'en\'efici\'e de l'aide vigilante de Myriam D\'echamps pour d\'ebusquer les failles et les erreurs. Je lui dois beaucoup, \`a la fois comme inspiratrice et comme premi\`ere lectrice et correctrice de cet article.

\section{D\'efinitions et principaux r\'esultats}

Soit $G$ un groupe ab\'elien compact et $\Gamma$ son dual, qui est un groupe ab\'elien discret. Soit $\Lambda$ une partie de $\Gamma$, et $S\geq 1$. On dit que $\Lambda$ est $S$--Sidon si, pour tout polyn\^ome ``trigonom\'etrique''
\begin{equation*}
P(g) = \sum_{\gamma\in \Lambda} a_\gamma \ \gamma(g)\qquad (a_\gamma \in \mathbb{C})\,,
\end{equation*}
on a
\begin{equation*}
\sum |a_\gamma| \leq S\ \sup_{g\in b} |P(g)|\,.
\end{equation*}
On dit que $\Lambda$ est Sidon s'il est $S$--Sidon pour un $S$  convenable. On conna\^{i}t maintenant un grand nombre de d\'efinitions \'equivalentes \cite{rudin} \cite{lr} \cite {lique}.

Parmi les ensembles de Sidon se trouvent les ensembles quasi--ind\'ependants, dont voici la d\'efinition : $\Lambda$ est quasi--ind\'ependant si la relation
\begin{equation*}
\sum_{\gamma\in \Lambda} \varepsilon_\gamma \ \gamma=0 \qquad (\varepsilon_\gamma \in \{-1,0,1\})
\end{equation*}
n'a lieu que lorsque tous les $\varepsilon_\gamma$ sont nuls.

On appellera maille dans $\Gamma$ tout ensemble de la forme
$$
M=M((\gamma_j)_{j=1,2,\ldots k}, E) = \sum_{j=1}^k \{n_j \gamma_j\}\,, \leqno{(1.1)}
$$
o\`u les $\gamma_j$ $(j=1,2,\ldots k)$ sont des \'el\'ements de $\Gamma$ et les $(n_j)_{j\in (1,2,\ldots k)}$ appartiennent \`a un ensemble $E$ dans $\mathbb{Z}^k$. Quand $M$ est de la forme (1.1) avec $|n_j| \leq h$ pour tout $j$, nous dirons que $M$ est une $k$--maille de hauteur $h$. Il y a g\'en\'eralement plusieurs \'ecritures de la forme $(1.1)$ pour un ensemble $M$ donn\'e ; $k$ et $h$ d\'ependent de l'\'ecriture.

Au sens g\'en\'eral, on dira qu'une partie $\Lambda$ de $\Gamma$ v\'erifie une condition de maille si l'on~a
$$
|\Lambda \cap M| \leq H(k,E) \leqno{(1.2)}
$$
pour toute maille $M$, $H(\cdot,\cdot)$ \'etant une fonction convenable. On a \'ecrit $|\Lambda\cap M|$ pour le cardinal de $\Lambda\cap M$, et ce sera la notation utilis\'ee dans la suite.

Si $\Lambda$ est Sidon, $\Lambda$ v\'erifie une condition de maille avec
$$
H(k,E) =Ck \log(1+\sup_{(n_j)\in E} (|n_1|+\cdots+|n_k|))\leqno{(1.3)}
$$
o\`u $C$ ne d\'epend que de $\Lambda$. De plus, si $\Lambda$ est $S$--Sidon, $C$ ne d\'epend que de $S$. C'est la condition de maille de~\cite{kahane}.

Nous allons voir dans la partie 2 que cette condition de maille est inam\'eliorable en plusieurs sens : on ne peut pas remplacer dans (1.3) la fonction $\log$ par une fonction qui soit $o(\log)$, ni remplacer la norme $\ell^1$ de $(n_1,n_2,\ldots n_k)$ par une norme substantiellement plus petite. Voici le r\'esultat.

\begin{montheo}
Soit $\Gamma$ un groupe ab\'elien discret contenant des \'el\'ements d'ordre arbitrairement grand.  Alors $\Gamma$ contient un ensemble quasi--ind\'ependant $\Lambda$
tel que pour tout entier $k$ il existe une $k$--maille $M$ de hauteur $1$ v\'erifiant
$$
|\Lambda \cap M| \geq {1\over4} k \log_2 k\,. \leqno{(1.4)}
$$
\end{montheo}

Ce th\'eor\`eme doit \^etre mis en rapport avec un r\'esultat de Pisier, la proposition 7.3 de \cite{pisier}, qui \'etablit une propri\'et\'e analogue lorsque $\Gamma$ est le groupe dual de $\mathbb{T}^N$ de fa\c con qualitative et non constructive, tandis que la preuve du th\'eor\`eme~1 est une construction explicite et \'el\'ementaire.

L'hypoth\`ese que $\Gamma$ contient des \'el\'ements d'ordre arbitrairement grand est essentielle. En effet, on sait par une autre proposition de Pisier que, dans le groupe $\Gamma$ dual de $G=\prod\limits_{j=1}^\infty(\mathbb{Z}/p_j\mathbb{Z})$, o\`u $(p_j)$ est une suite born\'ee d'entiers, on peut attacher  \`a tout ensemble de Sidon $\Lambda$ une constante $K$ telle que, pour tout sous--groupe fini $H$ de $\Gamma$, on ait
$$
|\Lambda\cap H| \leq K \ \textrm{rang}(H) \leqno{(1.5)}
$$
(corollaire 3.3 de \cite{pisier} ; quand $p_j = p$ fix\'e, c'est un r\'esultat de Malliavin--Malliavin \cite{mama}). Il s'ensuit que pour toute $k$--maille $M$ on a $|\Gamma\cap M| \le K k$.

Par ailleurs, la condition de maille est loin d'\^etre suffisante pour qu'un ensemble soit Sidon. M\^eme consid\'erablement renforc\'ee, elle ne garantit rien de tel. C'est ce que montrent les th\'eor\`emes 2 et 3, dont les preuves sont donn\'ees dans les parties 4 et 5.

\begin{montheo}
Soit $p$ premier, $\Gamma$ le groupe dual de $G=(\Z/p\Z)^\N$, et $w(x)$ une fonction croissante de $x(\geq 1)$, telle que
$$
w(x)\ge 1\qquad \hbox{et} \qquad \lim_{x\rg \infty} w(x) = \infty\,. \leqno{(1.6)}
$$
Il existe alors une partie $\Lambda$ de $\Gamma$, non Sidon, telle que pour tout entier $k\ge 1$ et toute $k$--maille $M$ on ait
$$
|\Lambda\cap M| \le k w(k)\,. \leqno{(1.7)}
$$
\end{montheo}

\begin{montheo}
Soit  $\Gamma=\Z$, et $w(x)$ comme dans le th\'eor\`eme~$2$. Il existe alors une partie $\Lambda$ de $\Z$, non Sidon, telle que, pour tout couple $(h,k)$ d'entiers positifs et toute $k$--maille $M$ de hauteur $h$ on ait
$$
|\Lambda\cap M| \le kw (kh) \leqno(1.8)
$$
\end{montheo}

La m\'ethode  de s\'election al\'eatoire utilis\'ee pour ces th\'eor\`emes et exprim\'ee par le lemme de la partie 3 est inspir\'ee de la note de Katznelson et Malliavin \cite{katmal} relative \`a la ``conjecture de dichotomie'' : ou bien $\Lambda$ est Sidon, ou bien c'est un ensemble d'analyticit\'e. La partie 6 rappelle la d\'efinition d'un ensemble d'analyticit\'e, et am\'eliore le th\'eor\`eme~2 dans le cas $p=2$ sous la forme que voici :

\begin{montheo}
Quand $p=2$, l'\'enonc\'e du th\'eor\`eme $2$ est valable en rempla\c cant ``non Sidon'' par ``d'analyticit\'e''.
\end{montheo}

Le cas $p=2$ n'a rien de sp\'ecial, sinon la relative facilit\'e d'\'ecriture des calculs.

L'appendice donne des estimations de distributions classiques, utilis\'ees dans l'article.

\section{Une construction d'ensembles quasi-ind\'e\-pendants. Preuve du th\'eor\`eme~1}

Nous allons d'abord nous placer dans la maille $\{-1,0,1\}^n$ de $\Z^n$ ($n$--maille de hauteur 1) et y construire un ensemble quasi--ind\'ependant $(q\cdot i\cdot)$ lorsque $n$ est une puissance de~2.

Lorsque $n=2$, les vecteurs colonnes de la matrice
$$
\left(\begin{matrix}
1 &1 &1\\
1 &-1 &0
\end{matrix} \right)
$$
sont $q\cdot i\cdot$ . V\'erifions--le en d\'etail. En effet, si
$$
\varepsilon_1 \left(\begin{matrix}
1\\
1
\end{matrix}\right) + 
\varepsilon_2 \left(\begin{matrix}
1\\
-1
\end{matrix}\right) +
\varepsilon_3 \left(\begin{matrix}
1\\
0
\end{matrix}\right)  =
 \left(\begin{matrix}
0\\
0
\end{matrix}\right) 
$$
avec  $\varepsilon_j\in \{-1,0,1\}$, on a d'abord $\varepsilon_1=\varepsilon_2$ (seconde ligne), puis $\varepsilon_3=0$ modulo 2 donc $\varepsilon_3=0$ (premi\`ere ligne), puis $\varepsilon_1=\varepsilon_2=0$ (ind\'ependance de $ \left(\begin{matrix}
1\\
1
\end{matrix}\right) $ et $ \left(\begin{matrix}
1\\
-1
\end{matrix}\right) $).

Lorsque $n=2^\nu$, on va construire par r\'ecurrence des matrices $A_\nu$ \`a $2^\nu$ lignes et $N_\nu$ colonnes, dont les colonnes sont dans $\{-1,0,1\}^n$ et sont $q\cdot i\cdot$. Pour $\nu=1$, c'est fait, avec $N_1=3$. On passe de $A_\nu$ \`a $A_{\nu+1}$ par le proc\'ed\'e figur\'e
$$
A_{\nu+1}=
\begin{tabular}{|c|c|c|}
\hline
$A_\nu$ & $A_\nu$ &$I_\nu$\\
\hline
$A_\nu $&$-A_\nu$ &0\\
\hline
\end{tabular}
$$
o\`u $I_\nu$ est la matrice unit\'e $2^\nu\times 2^\nu$. Montrons que les colonnes de $A_{\nu+1}$ sont $q\cdot i\cdot$ lorsque celles de $A_\nu$ le sont. Une relation lin\'eaire \`a coefficients $-1,0$ ou $1$ entre les colonnes de $A_{\nu+1}$ s'\'ecrit, en posant $n=2^\nu$ et $N=N_\nu$,
$$
\begin{matrix}
&A_\nu (\varepsilon_1^1, \varepsilon_1^2,\ldots \varepsilon_1^N)^t +A_\nu (\varepsilon_2^1, \varepsilon_2^2,\ldots \varepsilon_2^N)^t +I_\nu(\varepsilon_3^1,\varepsilon_3^2,\ldots\varepsilon_3^n)^t =0 \\
\noalign{\vskip2mm}
&A_\nu (\varepsilon_1^1, \varepsilon_1^2,\ldots \varepsilon_1^N)^t - 
A_\nu (\varepsilon_2^1, \varepsilon_2^2,\ldots \varepsilon_2^N)^t = 0\hfill
\end{matrix}
$$
$(\varepsilon_j^k \in \{-1,0,1\})$. En ajoutant, on voit que les lignes de $I_\nu(\varepsilon_3^1,\varepsilon_3^2,\ldots \varepsilon_3^n)^t$ sont nulles modulo 2, donc nulles, donc $\varepsilon_3^1=\varepsilon_3^2=\cdots=\varepsilon_3^n=0$. Il en r\'esulte
$$
A_\nu(\varepsilon_1^1,\varepsilon_1^2,\ldots \varepsilon_1^N)^t = A_\nu(\varepsilon_2^1,\varepsilon_2^2,\ldots \varepsilon_2^N)^t=0
$$
et la quasi--ind\'ependance des colonnes de $A_\nu$  entra\^{i}ne que tous les $\varepsilon_j^k$ sont nuls . Les colonnes de $A_{\nu+1}$ sont donc bien $q\cdot i\cdot$.

Calculons $N_\nu$. Partons de $N_0=1$. La construction donne
$$
N_\nu=2N_{\nu-1} +2^{\nu-1}
$$
soit
$$
2^{-\nu} N_\nu = 2^{-(\nu-1)} N_{\nu-1} + {1\over2} = \cdots = N_0 +{\nu\over2}
$$
donc
$$
N_\nu = 2^{\nu-1}(2+\nu)\,.
$$
Cela suffit \`a montrer que la condition de maille (1.2)--(1.3) est inam\'eliorable au sens pr\'ecis\'e dans la partie~1, et pour donner une minoration de la constante $C$ de (1.3) lorsque $\Lambda$ est $q\cdot i\cdot$ :
$$
c \ge {1\over2\ \log 2}\,.
$$

Pour d\'emontrer le th\'eor\`eme 1, on choisit dans $\Gamma$ une suite $(\beta_j)_{j\ge 1}$ tr\`es dissoci\'ee dans le sens suivant : il n'y a pas de relation lin\'eaire non triviale du type $\sum\, n_j\beta_j=0$ (somme finie) avec $n_j\in \Z$ et $|n_j| \le N_\nu$ quand $2^\nu \le j < 2^{\nu+1}$ $(\nu\ge 1)$. L'hypoth\`ese faite sur $\Gamma$ permet de construire une telle suite par r\'ecurrence. Pour
$$
\sum_{i=1}^{\nu-1} N_i < \ell \le \sum_{i=1}^\nu N_i
$$
on d\'efinit le vecteur ligne $(\gamma_\ell)$ comme
$$
(\gamma_\ell) = (\beta_{2^\nu}, \beta_{2^\nu+1},\ldots \beta_{2^{\nu+1}-1}) A_\nu \,.
$$
La suite cherch\'ee est $(\gamma_\ell)_{\ell\ge1}$. Elle est $q\cdot i\cdot$ parce que toute expression de la forme $\sum \varepsilon_\ell\gamma_\ell$ $(\varepsilon_\ell \in \{-1,0,1\})$ s'\'ecrit $\sum n_j \beta_j$ avec $|n_j| \le N_\nu$ quand $2^\nu \le j < 2^{\nu+1}$. Elle a $N_\nu$ termes dans la maille
$$
M= \{\sum \varepsilon_j \beta_j\,;\  \varepsilon_j \in \{-1,0,1\}\,;\ 2^\nu \le j < 2^{\nu+1}\}\,;
$$
$M$ est une $k$--maille de hauteur 1 lorsque $k\in [2^\nu, 2^{\nu+1}[$, et alors\break $N_\nu > {1\over4}\ k\ \log_2\ k$, ce qui \'etablit (1.4) et d\'emontre le th\'eor\`eme. \hfill $\blacksquare$

Si l'on se restreint aux valeurs de $k$ qui sont des puissances de 2, on peut minorer $|\Lambda \cap M|$ par ${1\over2}\ k\ \log_2\ k$ au lieu de ${1\over4}\ k\ \log_2\ k$.

\section{S\'elections et ind\'ependance dans $(\Z/p\Z)^\nu$. Un lemme}

Soit $p$ un nombre premier, $\nu$ un entier $\ge 1$, $X=(\Z/p\Z)^\nu$ et $(\Omega,P)$ un espace de probabilit\'e. Donnons--nous $\alpha$, $0<\alpha<1$, et associons--lui l'\'echantillon (= suite de $v\cdot a\cdot i\cdot i\cdot d\cdot)$ $\alpha(x,\omega)$ $(x\in X,\ \omega\in \Omega)$ de loi de Bernoulli $B(1,\alpha)$ et l'ensemble
$$
\Lambda(\omega) = \{x\in X : \alpha(x,\omega)=1\}\,. \leqno(3.1)
$$
Ainsi $|\Lambda (\omega)|$ a pour loi $B(p^\nu,\alpha)$ et l'on a (voir (7.9))
$$
P({1\over2}p ^\nu \alpha \le |\Lambda(\omega)| \le {3\over2} p^\nu \alpha) > 1- 2e^{-{1\over32}p^\nu\alpha} \leqno(3.2)
$$
Choisissons un entier $\ell$, $1\le \ell \le {p^\nu\over2\nu}$, et prenons
$$
\alpha=2\ell \nu p^{-\nu}\,; \leqno(3.3)
$$
ainsi
$$
P(\ell \nu \le |\Lambda(\omega)| \le 3\ell \nu) > 1-2 e^{-{1\over16}\ell\nu}\,. \leqno(3.4)
$$

Donnons--nous maintenant $\beta$, $0<\beta<1$, et associons--lui l'\'echantillon $\beta(x,\omega)$ $(x\in X, \omega\in \Omega)$ de loi $B(1,\beta)$, ind\'ependant des $\alpha(x,\omega)$. Soit 
$$
\lambda(\omega) = \{ x\in X : \alpha (x,\omega) \beta(x,\omega) = 1 \}\,. \leqno(3.5)
$$
Ainsi $\lambda(\omega)$ est une partie de $\Lambda(\omega)$ et $|\lambda(\omega)|$ a pour loi $B(p^\nu,\alpha\beta)$. Nous allons montrer que, si $\beta$ est bien choisi, la probabilit\'e que $\lambda(\omega)$ soit un syst\`eme libre dans l'espace vectoriel $X$ est voisine de~1.

Pour construire $\lambda(\omega)$, on peut commencer par choisir $k=k(\omega)$ al\'eatoire de loi $B(p^\nu,\alpha\beta)$, puis disposer au hasard sur $X$ $k$ points $\lambda_1,\lambda_2,\ldots \lambda_k$. Fixons $|\lambda|=k$ ; alors, pour $j<k$,
$$
P(\lambda_1,\ldots \lambda_{j+1} \ \hbox{\rm li\'e}\ |\lambda_1,\ldots\lambda_j\ \hbox{libre}) = {p^j-j\over p^\nu-j} < p^{j-\nu}
$$
donc
$$
P(\lambda\ \hbox{li\'e} \mid |\lambda|=k) <p^{-\nu}(1+p+\cdots +p^{k-1}) < p^{k-\nu}\,.
$$
Il s'ensuit que
$$
P(\lambda\ \hbox{li\'e}) < \sum_k P(|\lambda|=k) p^{k-\nu}
$$
et le second membre peut s'\'ecrire $p^{-\nu}(E(p^Z))^{p^\nu}$, $Z$ \'etant une $v\cdot a\cdot$ de loi $B(1,\alpha\beta)$, donc
$$
P(\lambda\ \hbox{li\'e}) < p^\nu (1-\alpha\beta+\alpha\beta p)^{p^\nu} < p^{-\nu}\exp (\alpha\beta(p-1)p^\nu)\,.
$$
Ce dernier terme est inf\'erieur \`a $p^{-\nu/2}$ si $\alpha\beta(p-1)p^\nu <{1\over2}\nu \log p$, soit, compte tenu de (3.3), $\beta <{1\over4}{\log p\over(p-1)\ell}$. Choisissons d\'esormais
$$
\beta={1\over4 p \ell}\,.
\leqno(3.6)
$$
Ainsi $P(\lambda\ \hbox{li\'e})<p^{-\nu/2}$.

D\'ecomposons $(\Omega,P)$ en un produit $(\Omega_\alpha,P_\alpha) \times (\Omega_\beta,P_\beta)$, les $\alpha(x,\omega)$ \'etant d\'efinis sur $\Omega_\alpha$ et les $\beta(x,\omega)$ sur $\Omega_\beta$. On~a
$$
E_\alpha(P_\beta (\lambda \ \hbox{li\'e}))= P(\lambda\ \hbox{li\'e}) <p^{-\nu/2}
$$
donc
$$
P_\alpha(P_\beta(\lambda \ \hbox{li\'e}) > p^{-\nu/4}) < p^{-\nu/4}
$$
et on sait, par (3.4), que
$$
P_\alpha(\ell \nu \le |\Lambda(\omega)| \le 3\ell \nu) > 1-2 e^{-{1\over16}\ell \nu}\,.
$$
Or
$$
p^{-\nu/4} < 1-2 e^{-{1\over16}\ell\nu}
$$
d\`es que $\nu\ge 16$. Sous cette condition, on peut choisir un point dans $\Omega_\alpha$, donc choisir $\Lambda$, de fa\c con que l'on ait \`a la fois $\ell\nu\le |\Lambda| \le 3\ell \nu$ et $P_\beta(\lambda \ \hbox{li\'e}) < p^{-\nu/4}$.

$\Lambda$ \'etant ainsi choisi, soit $A$ une partie de $\Lambda$. On a
$$
\begin{matrix}
P_\beta(A \subset \lambda) = \beta^{|A|}\,,\hfill \\
\noalign{\vskip2mm}
P_\beta(A\not\subset \lambda)\ \hbox{ou}\ \lambda\ \hbox{li\'e}) < 1-\beta^{|A|}+p^{-\nu/4}\,,
\end{matrix}
$$
donc, si $|A| \le K\nu$ et $\beta^K \ge p^{-1/4}$, $A$ est libre. La seconde condition est v\'erifi\'ee lorsque $K\le K_\ell$ avec (suivant (3.6))
$$
K_\ell = {1\over4} {\log p\over \log(4p\ell)}\leqno(3.7)
$$
Exprimons le r\'esultat.

\vskip2mm
{\bf Lemme}.
\textit{Soit $p$ premier, $\nu$ entier $\ge 16$ et $1\le \ell\le {p^\nu\over2\nu}$.  Il existe alors dans $(\Z/p\Z)^\nu$ une partie $\Lambda$ telle que $\ell\nu \le |\Lambda| \le 3\ell \nu$ et que toute partie de $\Lambda$ de cardinal inf\'erieur ou \'egal \`a $K_\ell \nu$ soit libre dans $(\Z/p\Z)^\nu$ (espace vectoriel sur $\Z/p\Z)$.}

\section{Preuve du th\'eor\`eme 2}

Ici $G=(\Z/p\Z)^\N$, $\Gamma$ est le dual de $G$, $G$ et $\Gamma$ sont des espaces vectoriels sur le corps $\Z/p\Z$. Soit $(\beta_i)\ (i\in \N)$ la base canonique de $\Gamma$, c'est--\`a--dire $\beta_i(x)=x_i$ quand $x=(x_0,x_1\cdots)\in G$. On r\'epartit les $\beta_i$ en blocs $B_\ell$ disjoints de cardinaux $|B_\ell|= \nu_\ell\ge 16$ tendant vers l'infini $(\ell = 2,3,4,\ldots)$. J'indiquerai plus loin (formule (4.3)) comment choisir les $\nu_\ell$ en fonction de $w(\cdot)$. Les $\beta_i\in B_\ell$ engendrent un sous--espace $\Gamma_\ell $ de $\Gamma$ isomorphe \`a $(\Z/p\Z)^{\nu_\ell}$, et d'apr\`es le lemme chaque $\Gamma_\ell$ contient un $\Lambda_\ell$ tel~que
$$
\ell \nu_\ell \le |\Lambda_\ell| \le 3 \ell \nu_\ell
$$
et que toute partie de $\Lambda_\ell$ de cardinal inf\'erieur ou \'egal \`a $K_\ell \nu_\ell$ est libre ($K_\ell$ \'etant d\'efini en~(3.7)).

D\'efinissons $\Lambda$ comme la r\'eunion des $\Lambda_\ell$. Comme
$$
|\Gamma_\ell \cap \Lambda| /\hbox{rang}\ \Gamma_\ell \ge \ell\,,
$$
qui n'est pas born\'e, $\Lambda$ n'est pas de Sidon (\cite{pisier}, \cite{mama}, voir (1.5)).

Soit $M$ une $k$--maille, $M\cap \Lambda=A$ et $M\cap \Lambda_\ell=A_\ell$. On veut montrer que
$$
|A| = \sum |A_\ell| \le k w (k)\,.
\leqno(4.1)
$$

Comme les $A_\ell$ appartiennent \`a des sous--espaces $\Gamma_\ell$ ind\'ependants, le rang de leur r\'eunion $A$ est la somme de leurs rangs. R\'epartissons les $\ell$ en deux  classes, $U$ et $W$, suivant que
$$
A_\ell < K_\ell \nu_\ell \qquad(\ell\in U)   
$$
ou
$$
  A_\ell \ge K_\ell \nu_\ell   \qquad(\ell \in W)\,.   
$$
Si $\ell\in U$ on a $\hbox{rang } A_\ell = |A_\ell|$ et, si $\ell\in W$, $\hbox{rang } A_\ell \ge K_\ell \nu_\ell$ puisque toute partie de $A_\ell$ de cardinal $\le K_\ell \nu_\ell$ est libre. Comme $\hbox{rang } A \le \hbox{rang } M \le k$, on~a
$$
k \ge \sum_U |A_\ell| + \sum_W K_\ell \nu_\ell\,.
$$
Si $W$ est vide, on a $|A|\le k$ et (4.1) est v\'erifi\'ee. Supposons donc $W$ non vide. On a toujours $|A_\ell| \le |\Lambda_\ell| \le 3\ell \nu_\ell$, donc
$$
\sum_W |A_\ell| \le \sum_W K_\ell \nu_\ell \sup_{W} {3\ell\over K_\ell}\,.
$$
et finalement
$$
|A| = \sum_U |A_\ell | + \sum_W |A_\ellÊ| \ \le k \sup_{W} {3\ell\over K_\ell}\,.
\leqno(4.2)
$$
Quand $\ell\in W$ on a
$$
k\ge \hbox{rang }A_\ell \ge K_\ell \nu_\ell
$$
donc (4.2) entra\^{i}ne (4.1) lorsque
$$
{3\ell\over K_\ell} \le w (K_\ell \nu_\ell)\,.
\leqno(4.3)
$$
C'est la condition que nous imposons pour le choix des $\nu_\ell$ ; elle garantit (4.1), ce qui ach\`eve la preuve du th\'eor\`eme~2. \hfill $\blacksquare$

\section{Preuve du th\'eor\`eme 3}

La preuve du th\'eor\`eme 3 s'inspire de celle du th\'eor\`eme~2, mais comme ici $\Gamma=\Z$ la notion de rang d'une partie de $\Gamma$ doit \^etre remplac\'ee par un substitut. La m\'ethode est \'evidente : elle consiste \`a \'etaler dans $\Z$ des copies d'ensembles du type $(\Z/p\Z)^\nu$ et d'y s\'electionner des parties de cardinal comparable \`a $\ell \nu$, pour des valeurs diff\'erentes de $p$, $\nu$, $\ell$. Mais on est forc\'e de pr\'eter attention aux d\'etails.

Pour toute partie finie $B$ de $\Gamma$ et tout entier impair $q$, d\'esignons par $V_q(B)$ l'ensemble des combinaisons lin\'eaires d'\'el\'ements de $B$ \`a coefficient entiers tels que $q\ge 2\ \sup\, |\hbox{coefficients}| +1$. C'est une $|B|$--maille de hauteur ${1\over2}(q-1)$. On dira que $V_q(B)$ est bien \'etal\'e si toutes ces combinaisons lin\'eaires sont distinctes. On a alors
$$
|V_q(B) |= q^{|B|}\,.
$$

On va d\'efinir en fonction de $w(\cdot)$ des suites croissantes au sens large $\ell_j\ (\ell_j >1)$, $p_j$ (nombres premiers) et $\nu_j$ (entiers $\ge 16$) $(j=1,2,\ldots)$. Observons que dans le lemme on peut remplacer $K_\ell$, donn\'e par (3.7), par ${1\over8}$ lorsque $4\ell<p$. Pour profiter de cette commodit\'e imposons
$$
4\ell_j <p_j
\leqno(5.1)
$$

Pour $\nu_1+\nu_2+\cdots+ \nu_{j-1} < i Ê\le \nu_1+\nu_2 +\cdots +\nu_j$ posons
$$
q(i) = 2\nu_j \left({p_j -1\over2}\right)^2 +1\,.
\leqno(5.2)
$$
D\'efinissons par induction une suite d'\'el\'ements $\beta_i$ de $\Gamma$ $(i=1,2,\ldots)$ assez rapidement croissante pour que les combinaisons lin\'eaires
$$
\sum m_i\beta_i,\ m_i\in \Z,\ |m_i| \le {1\over2}(q(i)-1)
\leqno(5.3)
$$
$(i=1,2,\ldots)$ soient toutes distinctes. Soit
$$
B_j = \{ \beta_i : \nu_1+\nu_2 + \cdots +\nu_{j-1} < i \le \nu_1+\nu_2 +\cdots \nu_j\}\,.
$$
Les conditions (5.2) et (5.3) impliquent que chaque $V_{p_j}(B_j)$ est bien \'etal\'e.

Appliquons la base canonique de $(\Z/p_j\Z)^{\nu_j}$ sur $B_j$, et $(\Z/p_j\Z)^{\nu_j}$, identifi\'e \`a l'ensemble des combinaisons lin\'eaires des \'el\'ements de la base canonique \`a coefficients entiers compris entre $-{1\over2}(p_j-1)$ et ${1\over2}(p_j-1)$, sur $V_{p_j} (B_j)$. D\'esignons par $\Lambda_j$ la partie de $V_{p_j}(B_j)$ qui est l'image dans cette application de la partie de $(\Z/p_j\Z)^{\nu_j}$ fournie par le lemme. Ainsi
$$
\ell_j \nu_j \le |\Lambda_j| \le 3\ell_j\nu_j
$$
et (en tenant compte de (5.1)) toute partie de $\Lambda_j$ de $\hbox{cardinal}\le {1\over8}\nu_j$ est l'image d'une partie libre de $(\Z/p_j\Z)^{\nu_j}$. Nous conviendrons d'\'ecrire qu'une telle image est ``ind\'ependante''.

V\'erifions que, si $A'$ est une partie ``ind\'ependante'' de $V_{p_j}(B_j)$ et si $p_j\ge p$ impair, on~a
$$
|V_p(A')| = p^{|A'|}\,.
\leqno(5.4)
$$
En effet, les \'el\'ements de $V_p(A')$ s'\'ecrivent $\sum\limits_{a\in A'}m_aa$, soit
$$
\sum_{a\in A'}m_a \sum_{\nu_1+\cdots+\nu_{j-1}<i\le \nu_1+\cdots+\nu_j} n_i(a)\beta_i\,,
\leqno(5.5)
$$
avec $|A'|\le \nu_j$, $|m_a| \le {p-1\over2}$ et $|n_i(a)| \le {p_j-1\over2}$, donc ils sont de la forme (5.3), et la construction des $\beta_i$ garantit que $V_p(A')$ est bien \'etal\'e. De plus, l'ind\'ependance assure que les combinaisons lin\'eaires dans $(\Z/p_j\Z)^{\nu_j}$ dont les (5.5) sont les images sont distinctes, d'o\`u~(5.4).

On aura besoin de (5.4) sous l'hypoth\`ese plus large que $A'$ est une r\'eunion finie de $A_j'$ qui sont des parties ``ind\'ependantes'' de $V_{p_j}(B_j)$ $(j=j_0,j_0+1,\ldots)$, avec $p_{j_0}\ge p$. En effet, les \'el\'ements de $V_p(A')$ sont toujours de la forme (5.3), et leur nombre est
$$
|V_p(A')| = \prod_j |V_p(A_j')|=p^{\sum |A_j'|}= p^{|A'|}\,.
$$

Posons d\'esormais $\Lambda=\bigcup\limits_{j\ge 1}\Lambda_j$. Etant donn\'e $M$, $k$--maille de hauteur $h$, soit
$$
A=\Lambda \cap M,\qquad A_j = \Lambda_j \cap M\,.
$$
Nous nous proposons de montrer que, par un choix convenable des suites $(p_j)$, $(\nu_j)$ et $(\ell_j)$, ne d\'ependant que de $w(\cdot)$, on a la conclusion du th\'eor\`eme~3, c'est--\`a--dire
$$
|A| = \sum |A_j| \le k w (kh)\,.
\leqno(5.6)
$$
Il nous restera \`a v\'erifier ensuite que $\Lambda$ n'est pas Sidon.

Estimons s\'epar\'ement les sommes des $A_j$ correspondant \`a $j\le k$ et \`a $j>k$.

Comme $|A_j| \le |\Lambda_j| \le 3 \ell_j \nu_j$, on a
$$
\sum_{j\le k} |A_j| \le 3\ell_k \sum_{j\le k} \nu_j\,.
\leqno(5.7)
$$
Pour $j>k$, d\'esignons par $A_j'$ une partie de $A_j$ ``ind\'ependante'' maximale. Distinguons les deux cas :
$$
\begin{array}{ccc}
  &U : |A_j| \le {1\over8} \nu_j  \hfill \\
  \noalign{\vskip2mm}
  &W : {1\over8} \nu_j < |A_j| \le 3\nu_j\ell_j   \\
\end{array}
$$
D'apr\`es le lemme, $A_j'=A_j$ dans le cas $U$ et $|A_j'| \ge {1\over8}\nu_j$ dans le cas $W$ (gr\^ace \`a (5.1)). D\'ecomposons en cons\'equence la somme $\sum\limits_{j>k}$ en $\sum\limits_U + \sum\limits_W$ :
$$
\begin{array}{lll}
  &\displaystyle\sum\limits_U |A_j| = \sum\limits_U |A_j'|  \\
  \noalign{\vskip2mm}
  &\displaystyle \sum_W |A_j| \le 3 \sum_W \nu_j\ell_j \le 24 \sum\limits_W \ell_j |A_j'| \le 24 \sup\limits_W \ell_j \sum\limits_W |A_j'| \\
\end{array}
$$
Comme $|A_j'| \le |M| \le (2h+1)^k$, on a dans le cas $W$ $\nu_j \le 8 (2h+1)^k$. Finalement
$$
\left\{
\begin{array}{lll}
  &\displaystyle \sum_{j>k} |A_j| \le X \sum_{j>k} |A_j'|   \\
  \noalign{\vskip2mm}
  &X =\displaystyle \sup (1,\ \sup_{\nu_j \le 8 (2h+1)^k} \ell_j  )\\ 
\end{array}
\right. \leqno(5.8)
$$

Pour utiliser (5.8), posons $A' = \bigcup\limits_{j>k} A_j'$ et $p=p_k$. La formule (5.4) s'applique : $|V_p(A')| =p^{|A'|}$. D'autre part les \'el\'ements de $V_p(A')$ s'\'ecrivent $\sum\limits_{a\in A'}m_a a$ avec $|m_a| \le {1\over2}(p-1)$ et, si la base de la maille $M$ est $(\gamma_1,\gamma_2,\ldots \gamma_k)$, chaque $a\in A'$ s'\'ecrit $\sum\limits_{1\le i\le k}n_i(a)\gamma_i$ avec $|n_i(a)|\le h$. Ainsi les \'el\'ements de $V_p(A')$ sont de la forme
$$
\sum_{1\le i\le k}\  \sum_{a\in A'} m_a n_i(a)\gamma_i\,.
$$
Les coefficients des $\gamma_i$ sont major\'es en module par ${1\over2}(p-1)h|A'|$. Leur nombre, pour un $i$ fix\'e, est major\'e par $ph|A'|$, donc
$$
|V_p(A')| \le (ph|A'|)^k\,.
\leqno(5.9)
$$
L'\'evaluation de $|A'|$ repose sur l'in\'egalit\'e, venant de (5.4) et (5.9),
$$
p^{|A'|} \le (ph |A'|)^k\,.
$$
En majorant $|A'|$ dans le second membre par $(2h+1)^k$, on obtient
$$
|A'| \le k \Big( 1+ {(k+1) \log (2h+1)\over\log p}\Big) \qquad (p=p_k)
\leqno(5.10)
$$

Reste \`a regrouper (5.7), (5.8) et (5.10) pour obtenir (5.6), moyennant un choix convenable des $p_j$, $\nu_j$ et $\ell_j$.

On choisit $p_k= p_k(h,k)$ de fa\c con que (5.10) entra\^{i}ne
$$
|A'| \le {1\over2} k\ w^{1/2}(hk)\,,
$$
les $\nu_j=\nu_j(h,k)$ de fa\c con que 
$$
\sum_{j\le k} \nu_j \le {1\over2} k\ w^{1/2} (hk)\,,
$$
et enfin les $\ell_j =\ell_j(h,k)$ de fa\c con que 1\up{o}) $4\ell_j < p_j$ (c'est la condition (5.1)) \ 2\up{o}) $3\ell_k \le w^{1/2}(hk)$, de fa\c con que (5.7) entra\^{i}ne
$$
\sum_{j\le k} |A_j| \le {1\over2} k \ w (hk)\,,
$$
3\up{o}) $X$ dans (5.8) v\'erifie $X\le w^{1/2}(hk)$, de fa\c con que
$$
\sum_{j>k} |A_j| \le {1\over2} k\ w(hk)\,.
$$
En gros, les $p_j$ croissent tr\`es vite, les $\nu_j$  lentement et les $\ell_j$ tr\`es lentement. On a obtenu
$$
\sum |A_j| \le k\ w(hk)\,,
$$
l'in\'egalit\'e (5.6) voulue.

Pour voir que $\Lambda$ n'est pas Sidon, il suffit de v\'erifier que le crit\`ere de Pisier (\cite{pisier}, \cite{lique} p.~483) n'est pas v\'erifi\'e, c'est--\`a--dire que pour tout $\delta >0$ il existe une partie de $\Lambda$, soit $\Lambda_\delta $, dont toute partie quasi--ind\'ependante a moins de $\delta |\Lambda_\delta |$ \'el\'ements. Cela a bien lieu en prenant pour $\Lambda_\delta$ 
 un $\Lambda_j$ avec $j$ assez grand \hfill$\blacksquare$
 
 \section{Ensembles d'analyticit\'e. Preuve du\break th\'eor\`eme~4}

On d\'esigne par $A(\Gamma)$ l'ensemble des transform\'ees de Fourier des fonctions int\'egrables sur $G : A(\Gamma)=FL^1(G)$. Quand $\Lambda\subset \Gamma$ on d\'esigne par $A(\Lambda)$ l'espace des restrictions \`a $\Lambda$ des fonctions appartenant \`a $A(\Gamma)$. On a toujours $A(\Lambda)\subset c_0(\Lambda)$, espace des fonctions d\'efinies sur $\Lambda$ et tendant vers $0$ \`a l'infini. L'une des d\'efinitions des ensembles de Sidon est l'\'egalit\'e de ces espaces comme ensembles : $A(\Lambda)=c_0(\Lambda)$.

On dit qu'une fonction $F$ d\'efinie sur un intervalle r\'eel ouvert $I$ contenant $O$ ``op\`ere dans $A(\Lambda)$'' si, pour toute $f\in A(\Lambda)$ \`a valeurs dans $I$ on a $F\circ f \in A(\Lambda)$. Il est n\'ecessaire pour cela que $F(O)=O$ et que $F$ soit continue en $O$, et c'est suffisant lorsque $\Lambda$ est Sidon. En tout cas il est suffisant que $F(O)=O$ et que $F$ soit analytique au voisinage de~$O$.

On dit que $\Lambda$ est ``ensemble d'analyticit\'e'' si les seules fonctions qui op\`erent dans $A(\Lambda)$ sont les fonctions analytiques nulles en $O$. La conjecture de dichotomie de Katznelson est qu'une partie de $\Gamma$ est soit Sidon, soit ensemble d'analyticit\'e (\cite{kakat} p. 112). L'article \cite{katmal} de Katznelson et Malliavin est une v\'erification de cette conjecture dans un cadre al\'eatoire sous la forme d'une propri\'et\'e presque s\^ure. Le th\'eor\`eme 4 est une variante de leurs r\'esultats. On s'est born\'e au cas $G=(\Z/2\Z)^\N$ non seulement pour simplifier les \'ecritures, mais aussi parce que c'est le cadre naturel pour tester la conjecture de Katznelson, toujours ouverte.

Revenons donc \`a la partie 3, avec maintenant $X=(\Z/2\Z)^\nu$ ; $(\Omega,P)$ est un espace de probabilit\'e, $\alpha(x,\omega)$ $(x\in X,\ \omega\in \Omega)$ un \'echantillon de loi $B(1,\alpha)$, et $\Lambda=\Lambda(\omega)=\{x\in X : \alpha(x,\omega)=1\}$.

Soit $Y$ le groupe des caract\`eres sur $X$, not\'e multiplicativement, et $\sigma=\sigma(\omega)$ la mesure de d\'ecompte sur~$\Lambda$ :
$$
\sigma=\sum_{x\in \Lambda}\delta _x\,.
$$
La transform\'ee de Fourier de $\sigma$ est 
$$
\hat\sigma(y,\omega) = \sum_{x\in \Lambda(\omega)} y(x) = \sum_{x\in X} \alpha(x,\omega)y(x)\,.
\leqno(6.1)
$$
Pour $y=1$, on trouve $|\Lambda(\omega)|$, et (3.2) s'\'ecrit
$$
P\Big( {1\over2} 2^\nu \alpha \le \hat \sigma (1,\omega) \le {3\over2} 2^\nu \alpha\Big) > 1-2 \exp (-2^{\nu-5}\alpha)\,.
\leqno(6.2)
$$
Pour $y\neq 1$, il y a autant de $x$ pour lesquels $y(x)=1$ que pour $y(x)=-1$. La derni\`ere somme dans (6.1) est donc la diff\'erence de deux $v\cdot a\cdot$ ind\'ependantes de lois $B(2^{\nu-1},\alpha)$. D'apr\`es (7.10) et (7.11),  si
$$
0<\lambda \le \big(2^{\nu-1}\alpha(1-\alpha)\big)^{1/2}
\leqno(6.3)
$$
on a
$$
P(|\hat\sigma(y,\omega)| > 2\lambda (2^\nu \alpha(1-\alpha))^{1/2}) \le 2\ e^{-{1\over2}\lambda^2}\,, 
$$
d'o\`u
$$
P\Big(\sup_{y\neq1}
|\hat\sigma(y,\omega)| > 2\lambda (2^\nu \alpha)^{1/2}\Big) \le 2^{\nu+1} e^{-{1\over2}\lambda^2}\,.
\leqno(6.4)
$$
Choisissons comme dans la partie 3 $\alpha=2\ell \nu 2^{-\nu}$ puis $K < {1\over4}{\log 2\over\log 8\ell}$ (cf. formules (3.3) et (3.7)), ce qui assure qu'avec une probabilit\'e sup\'erieure \`a une puissance n\'egative de $\nu$, \`a savoir $(8\ell)^{-K\nu}- 2^{-\nu/4}$, toute partie de $\Lambda(\omega)$ de cardinal$\ \le K\nu$ est libre. Choisissons $\lambda=10\nu^{1/2}$ (valeur permise par (6.3)). Alors, avec une probabilit\'e positive on~a
$$
\left\{
\begin{array}{ll}
  \hat\sigma(1,\omega) \ge \ell \nu   \\
  \noalign{\vskip1mm}
   \sup\limits_{y\neq1} |\hat\sigma (y,\omega)| \le 20 \nu^{1/2}(\ell \nu)^{1/2}   \\
\end{array}
\right.
$$
donc $\sigma$ est une mesure positive, port\'ee par $\Lambda$, telle que
$$
\sup_{y\neq1} |\hat\sigma(y,\omega)| \le {20\over\sqrt{\ell}} \ \hat\sigma(1,\omega)\,.
\leqno(6.5)
$$
On choisit $\ell>400$ et on fera ensuite tendre $\ell$ vers l'infini. On sait qu'une in\'egalit\'e du type (6.5) entra\^{i}ne que la r\'eunion des $\Lambda$ correspondants est un ensemble d'analyticit\'e \cite{kakat}. D\'etaillons le calcul, qui  est facile. D\'esormais $\omega$ est choisi pour avoir (6.5), et on \'ecrit $\sigma$ pour $\sigma(\omega)$.

Choisissons un entier $\rho<\nu$, et consid\'erons $\rho$ caract\`eres ind\'ependants $y_1, y_2,\ldots y_\rho$. Posons
$$
\begin{array}{lll}
 f = y_1 + y_2+\cdots + y_\rho   \\
v = \exp \big( i {\pi\over4}f) = 2^{-\rho/2}(1+y_1)(1+y_2)\cdots(1+y_\rho)  \\
\mu = vÊ\sigma\,.
\end{array}
$$
On a pour tout $y\in Y$
$$
\begin{array}{llll}
 \hat\mu(y) &= \displaystyle\sum_{y'y''=y} \hat v(y') \hat\sigma(y'')   \\
  &=\hat v(y) \hat\sigma(1)+\displaystyle\sum_{y'y''=y,\, y''\neq1} \hat v(y') \hat\sigma(y'')   \\
  &\le 2^{-\rho/2} \hat\sigma(1) +2^{\rho/2} \sup\limits_{y\neq1} |\hat\sigma(y)|\\
  &\le \hat\sigma(1) \big(2^{-\rho/2} + {20\over\sqrt{\ell}}\ 2^{\rho/2}\big)   
\end{array}
$$
Or, comme $\sigma=\mu v^{-1}$,
$$
\hat\sigma(1) \le \sup_{y} \ |\hat\mu(y)|\  \|v^{-1}\|_{A(\Lambda)}\,,
$$
d'o\`u
$$
\|v\|_{A(\Lambda)}= \|v^{-1}\|_{A(\Lambda)} \ge \Big(2^{-\rho/2} +{20\over\sqrt{\ell}} 2^{\rho/2}\Big)^{-1}
$$
Choisissons $\rho=\log_2{\sqrt{\ell}\over20}$ ; on obtient
$$
\|v\|_{A(\Lambda)} \ge {1\over2}\ 2^{\rho/2}\,.
$$
Comme $\|f\|_{A(\Lambda)} \le \| f\|_{A(X)}=\rho$, on a finalement
$$
\left\{
\begin{array}{ll}
  \|f\|_{A(\Lambda)} \le \rho  \\
  \noalign{\vskip2mm}
 \|\exp \displaystyle{i\pi\over4} f\|_{A(\Lambda)} \ge {1\over2} 2 ^{\rho/2}  \\  
\end{array}
\right.
\leqno(6.6)
$$

Reportons--nous maintenant \`a la partie 4, dans laquelle on a construit $\Lambda$ comme une r\'eunion de parties $\Lambda_\ell$ de blocs $\Gamma_\ell$ isomorphes \`a $(\Z/p\Z)^{\nu_\ell}$ $(\ell=2,3,4,\ldots)$ ; moyennant un choix convenable des $\nu_\ell$, \`a savoir (4.3), la conclusion du th\'eor\`eme~2 est valide. En prenant ci--dessus $\nu=\nu_\ell$ et en transportant la fonction $f$ sur $\Gamma_\ell$, on obtient une fonction $f_\ell\in A(\Lambda)$ pour laquelle (6.6) a lieu avec $\rho=\log_2 {\sqrt{\ell}\over20}$, et d'apr\`es un crit\`ere connu \cite{kakat} cela montre que $\Lambda$ est un ensemble d'analyticit\'e. \hfill$\blacksquare$

\section{Appendice : estimation de distributions}

On appelle sous--gaussienne de type $\tau$ une variable al\'eatoire centr\'ee $X$ dont la transform\'ee de Laplace v\'erifie
$$
E(e^{uX}) \le e^{{1\over2} u^2\tau^2} \qquad (u\in \R)\,.
$$
On sait que cela donne un bon contr\^ole de la distribution
$$
\left\{
\begin{array}{ll}
  P(X> \lambda\tau) &\le \inf\limits_{u} \exp \Big({1\over2} u^2 \tau^2 -\lambda u \tau\Big) = e^{-{1\over2}\lambda^2}   \\
   P(|X\ >\lambda \tau) &\le 2 e ^{-{1\over2}\lambda^2}   \\
\end{array}
\right.
\leqno(7.1)
$$
Il est int\'eressant pour certains calculs de disposer d'une version locale, que voici. On dira qu'une $v\cdot a\cdot$ centr\'ee $X$ est sous gaussienne de type $\tau/\!/ h$ (de type $\tau$ relativement \`a l'intervalle $(-h,h))$ si
$$
E(e^{uX}) \le e^{{1\over2} u^2\tau^2} \quad\hbox{quand }\ -h \le u\le h\,.
\leqno(7.2)
$$
Le calcul pr\'ec\'edent montre que les in\'egalit\'es (7.1) sont valables lorsque
$$
0 < \lambda < \tau h\,.
\leqno(7.3)
$$

Premi\`ere application, aux sommes de $v\cdot a \cdot$ ind\'ependantes.

Si les $v\cdot a \cdot$ centr\'ees $X_j$ sont sous--gaussiennes de types $\tau_j/\!/ h$ et ind\'ependantes, leur somme $X= X_1 +\cdots + X_N$ est sous--gaussienne de type $\tau/\!/ h$ avec $\tau^2 =\tau_1^2 +\cdots + \tau_N^2$, et (7.1) s'\'ecrit
$$
\left\{
\begin{array}{lll}
  P\Big(X_1 +\cdots + X_N    &> \lambda(\tau_1^2 +\cdots +\tau_N^2)^{1/2}\Big)   &\le e ^{-{1\over2}\lambda^2}   \\
P\Big(|X_1 +\cdots + X_N |   &> \lambda(\tau_1^2 +\cdots +\tau_N^2)^{1/2}\Big)   &\le 2e ^{-{1\over2}\lambda^2}    \\  
\end{array}
\right.
\leqno(7.4)
$$

Seconde application, aux  $v\cdot a \cdot$ de Bernoulli.

Soit $0 < \alpha<1$. Commen\c cons par v\'erifier l'in\'egalit\'e
$$
\alpha e^{(1-\alpha)u} + (1-\alpha) e^{-\alpha u} \le e^{2\alpha(1-\alpha)u^2}
\leqno(7.5)
$$
sous la condition
$$
- {1\over|2-4\alpha|} \le u \le {1\over |2-4\alpha|}\, .
\leqno(7.6)
$$

Quand $\alpha=0$ ou $1$, (7.5) a lieu, et en tous cas (7.5) s'\'ecrit
$$
\alpha e^u +1-\alpha \le \exp (2\alpha(1-\alpha) u^2 +\alpha u)\,.
$$
Sous cette forme, le premier membre est une fonction affine de $\alpha$, et il suffit de v\'erifier que le second membre est concave sur $[0,1]$ quand $u$ est fix\'e selon (7.6). Or ce second membre est de la forme $e^A$ et la condition de concavit\'e est $A'^2 +A'' \le 0$, soit ici
$$
(( 2-4\alpha) u^2 +u)^2 -4 u^2 \le 0
$$
ou encore
$$
((2-4\alpha) u+3) ((2-4\alpha)u-1)\le 0
$$
ce qui a lieu d'apr\`es (7.6).

Nous venons de montrer que, si $Y$ est une  $v\cdot a \cdot$ de loi $B(1,\alpha)$, $Y-\alpha$ est sous--gaussienne de type $2\sqrt{\alpha(1-\alpha)}/\!/ {1\over|2-4\alpha|}$.

Il r\'esulte de (7.4) que, si $Y$ est une  $v\cdot a \cdot$ de loi $B(N,\alpha)$, on a
$$
\left\{
\begin{array}{lll}
  P(Y-N\alpha &> 2\lambda \sqrt{N\alpha(1-\alpha)})   &\le e^{-{1\over2}\lambda^2}   \\
  \noalign{\vskip1mm}
  P(|Y-N\alpha| &> 2\lambda \sqrt{N\alpha(1-\alpha)})   &\le 2e^{-{1\over2}\lambda^2}   \\
\end{array}
\right.
\leqno(7.7)
$$
lorsque
$$
0<\lambda < {\sqrt{N\alpha(1-\alpha)}\over |1-2\alpha|}
\leqno(7.8)
$$
et en particulier quand $0<\lambda \le \sqrt{N\alpha(1-\alpha)}$.

Supposons $0<\alpha<{1\over2}$. Le choix de $\lambda={1\over4} \sqrt{N\alpha}$ donne
$$
P(|Y-N\alpha| >{1\over2} N\alpha) \le 2\ e^{-{1\over32} N\alpha\,.}
\leqno(7.9)
$$
Nous nous sommes servis de cette in\'egalit\'e dans la partie 3.

Consid\'erons enfin une $v\cdot a \cdot$ $Z$ de la forme $Y-Y'$, o\`u $Y$ et $Y'$ sont deux  $v\cdot a \cdot$ ind\'ependantes de m\^eme loi $B(N,\alpha)$. C'est une $v\cdot a \cdot$ centr\'ee, sous--gaussienne de type $2 \sqrt{2N\alpha(1-\alpha)}/\!/ {1\over|2-4\alpha|}$ $(0 < \alpha <1)$, donc
$$
P(|Z| > 2\lambda \sqrt{2N\alpha(1-\alpha)} < 2\ e^{-{1\over2} \lambda^2}
\leqno(7.10)
$$
sous la condition
$$
\lambda< \sqrt{{N\alpha(1-\alpha)\over |1-2\alpha|}}\,.
\leqno(7.11)
$$
Nous nous sommes servis de cela dans la partie 6.

\eject

\vskip4mm

\hfill\begin{minipage}{6,5cm}
Jean--Pierre Kahane

Laboratoire de Math\'ematique

Universit\'e Paris--Sud, B\^at. 425

91405 Orsay Cedex

\textsf{Jean-Pierre.Kahane@math.u-psud.fr}

31/08/2007

\end{minipage}

 \end{document}